\documentclass{article}

\usepackage[top=20truemm,bottom=20truemm,left=20truemm,right=20truemm]{geometry}
\usepackage{authblk}
\usepackage{abstract}
\usepackage{amsmath}
\usepackage{amssymb}
\usepackage{hyperref}
\usepackage[dvipdfmx]{graphicx}
\usepackage{color}

\author[1]{Takuya Tsuchiya\thanks{\href{mailto:t-tsuchiya@hi-tech.ac.jp}{\nolinkurl{t-tsuchiya@hi-tech.ac.jp}}}}
\author[2]{Makoto Nakamura}
\affil[1]{%
  Center for Liberal Arts and Sciences,
  Hachinohe Institute of Technology,
  Japan
}
\affil[2]{%
  Department of Pure and Applied Mathematics,
  Graduate School of Information Science and Technology,
  Osaka University,
  Japan
}

\title{Numerical accuracy and stability of semilinear Klein--Gordon equation
  in de Sitter spacetime}

\begin{document}

\maketitle

\begin{abstract}%
  Numerical simulations of the semilinear Klein--Gordon equation in the de
  Sitter spacetime are performed.
  We use two structure-preserving discrete forms of the Klein--Gordon equation.
  The disparity between the two forms is the discretization of the differential
  term.
  We show that one of the forms has higher numerical stability and second-order
  numerical accuracy with respect to the grid, and we explain the reason for the
  instability of the other form.
\end{abstract}

\section{Introduction}

The Klein--Gordon equation is one of the relativistic wave equations.
There have been some analytical investigations of this equation (e.g., \cite{Yagdjian-Galstian, Nakamura-2014-JMAA,Nakamura-2021-JMP})
.
However, it is difficult to 
quantitatively evaluate the solutions
analytically; therefore, we carry out numerical simulations to investigate the
solutions.
In this paper, we adopt the structure-preserving scheme \cite{
  Furihata-Matsuo-2010} as a discretized scheme to realize high-accuracy and
high-stability simulations.
Here, the word ``accuracy'' means the difference between the initial and
time-evolved values of constraints.
A value is, for example, the total Hamiltonian in the Hamiltonian system.
In general, numerical accuracy means the difference between the exact and
numerical solutions.
However, the exact solution is usually not obtained for the targets in numerical
calculations.
Thus, the conservation of the constraint is treated as the numerical accuracy
of the system in this paper.
In addition, the word ``stability'' means that the numerical solutions have no
numerical vibrations.
Precisely, the solution with numerical vibrations is less stable than that
without such vibrations.

In our previous paper \cite{Tsuchiya-Nakamura-2019-JCAM}, we reported some
numerical results of the solution of the Klein--Gordon equation with the
structure-preserving scheme.
However, we have recently reported in \cite{Tsuchiya-Nakamura-2022-ISAAC} that
numerical stability is not sufficient for the quantitative evaluations of the
numerical results in \cite{Tsuchiya-Nakamura-2019-JCAM}.
Thus, we propose another structure-preserving discretized form with higher
stability and comparable accuracy to the form in \cite{Tsuchiya-Nakamura-2019-JCAM}.

In this paper, we set the physical constants of the speed of light, the constant
of gravitation, and Dirac's constant as units.
Indices such as $(i, j, k, \dots)$ run from 1 to 3.
We use the Einstein convention of summation of repeated up--down indices.

\section{Semilinear Klein--Gordon equation in de Sitter spacetime}

The semilinear Klein--Gordon equation in the de Sitter spacetime (e.g.,
\cite{Tsuchiya-Nakamura-2019-JCAM}) is given by
\begin{align}
  -\partial_t^2\phi - nH \partial_t\phi + e^{-2Ht}\delta^{ij}(\partial_i
  \partial_j\phi) - m^2\phi = |\phi|^{p-1}\phi,
  \label{eq:KG1}
\end{align}
where $\phi$ is the dynamical variable, $n$ is the spacial dimension, $m$ is
the mass, $p$ is the integer greater than or equal to 2, and $H$ is the Hubble
constant.
$H>0$ means the expansion of space; conversely, $H<0$ means the contraction of
space.
$H=0$ means a flat space.

The Hamiltonian of \eqref{eq:KG1} is given by
\begin{align}
  \mathcal{H}
  &=
  \dfrac{1}{2}e^{-nHt}\psi^2
  + \dfrac{m^2}{2}e^{nHt}\phi^2
  + \dfrac{1}{p+1}e^{nHt}|\phi|^{p+1}
  + \dfrac{1}{2}e^{(n-2)Ht}\delta^{ij}(\partial_i\phi)(\partial_j\phi).
  \label{eq:HamiltonianDeSitter}
\end{align}
Then, the evolution equations are
\begin{align}
  \partial_t\phi
  &= e^{-nHt}\psi,
  \label{eq:phi}
  \\
  \partial_t\psi
  &=
  - m^2e^{nHt}\phi
  - e^{nHt}|\phi|^{p-1}\phi
  + e^{(n-2)Ht}\delta^{ij}(\partial_j\partial_i\phi).
  \label{eq:psi}
\end{align}
The total Hamiltonian $H_C$ is defined as
\begin{align}
  H_C(t):=
  \int_{\mathbb{R}^n} \mathcal{H}(t,x^i)\,d^{n}x,
  \label{eq:totalH}
\end{align}
and with \eqref{eq:phi}--\eqref{eq:psi}, the time derivative of
\eqref{eq:totalH} is
\begin{align}
  \partial_tH_C
  &=
  \dfrac{H}{2}
  e^{nHt}\int_{\mathbb{R}^n} d^n x\biggl(
  - ne^{-2nHt}\psi^2
  + nm^2\phi^2
  + \dfrac{2n}{p+1}|\phi|^{p+1}
  + (n-2)\delta^{ij}e^{-2Ht}(\partial_i\phi)(\partial_j\phi)
  \biggr)
  \nonumber\\
  &\quad
  + \int_{\mathbb{R}^n}d^nx\partial_j\left(
  c^2\delta^{ij}e^{-2Ht}\psi(\partial_i\phi)
  \right).
  \label{eq:timeDeriv_totalHC}
\end{align}
Note that $H$ is the Hubble constant and $H_C$ is the total Hamiltonian.
If we set an appropriate boundary condition such as the periodic boundary
condition, the last term on the right-hand side of \eqref{eq:timeDeriv_totalHC}
is eliminated.
In addition, if $H=0$, $\partial_tH_C=0$.
This means that $H_C$ is constant with time if $H=0$.
On the other hand, $H_C$ is not constant with time if $H\neq 0$.
To make a value constant with time, we define a value as
\begin{align}
  \tilde{H}_C(t)
  &:= H_C(t) - \int^t_0\partial_sH_C(s)ds
  \nonumber\\
  &=H_C(t)
  - \dfrac{H}{2}\int^t_0ds\,e^{nHs}
  \biggl\{
  \int_{\mathbb{R}^n} d^n x\biggl(
  - ne^{-2nHs}\psi^2
  + nm^2\phi^2
  + \dfrac{2n}{p+1}|\phi|^{p+1}
  + (n-2)e^{-2Hs}\delta^{ij}(\partial_i\phi)
  (\partial_j\phi)
  \biggr)
  \biggr\}.
  \label{eq:modTotalHam}
\end{align}
We call this value the modified total Hamiltonian, which is constant with time
if $H\neq 0$.
We judge whether the simulations are a success or a failure if $H_C$ in $H=0$
or $\tilde{H}_C$ in $H\neq 0$ is conserved or not in the time evolution.

\section{Discretized form of semilinear Klein--Gordon equation in de Sitter
  spacetime}

In this paper, we use two discretized forms.
The first form includes the product of the first-order central difference
formulae in the evolution equations.
The second form includes the product of the first-order forward and backward
difference formulae in the equations.
We call the first form Form I and the second one Form II.
Form I was suggested in \cite{Tsuchiya-Nakamura-2019-JCAM}, and Form II is a
new discretized form of the Klein--Gordon equation.

\subsection{Form I}

The discretized \eqref{eq:HamiltonianDeSitter}, \eqref{eq:phi}, and
\eqref{eq:psi} can be respectively defined as
\begin{align*}
  \mathcal{H}^{\mathrm{I}}{}^{(\ell)}_{(\boldsymbol{k})}
  &:= \frac{1}{2}e^{-nHt^{(\ell)}}
  (\psi^{(\ell)}_{(\boldsymbol{k})})^2
  + \dfrac{m^2}{2}e^{nHt^{(\ell)}}
  (\phi^{(\ell)}_{(\boldsymbol{k})})^2
  + \dfrac{1}{p+1}e^{nHt^{(\ell)}}
  |\phi^{(\ell)}_{(\boldsymbol{k})}|^{p+1}
  + \dfrac{1}{2}e^{(n-2)Ht^{(\ell)}}
  \delta^{ij}
  (\widehat{\delta}^{\langle1\rangle}_i
  \phi^{(\ell)}_{(\boldsymbol{k})})(\widehat{\delta}^{\langle1\rangle}_j
  \phi^{(\ell)}_{(\boldsymbol{k})}),
\end{align*}
\begin{align}
  \frac{\phi^{(\ell+1)}_{(\boldsymbol{k})}
    - \phi^{(\ell)}_{(\boldsymbol{k})}}{\Delta t}
  &:=
  \dfrac{1}{4}(e^{-nHt^{(\ell+1)}} + e^{-nHt^{(\ell)}})
  (\psi^{(\ell+1)}_{(\boldsymbol{k})} + \psi^{(\ell)}_{(\boldsymbol{k})}),
  \label{eq:DiscretePhi-I}
  \\
  \frac{\psi^{(\ell+1)}_{(\boldsymbol{k})}
    - \psi^{(\ell)}_{(\boldsymbol{k})}}{\Delta t}
  &:=
  - \frac{m^2}{4}(e^{nHt^{(\ell+1)}} + e^{nHt^{(\ell)}})
  (\phi^{(\ell+1)}_{(\boldsymbol{k})} + \phi^{(\ell)}_{(\boldsymbol{k})})
  - \frac{1}{2(p+1)}(e^{nHt^{(\ell+1)}}\!\! +\!e^{nHt^{(\ell)}})
  \frac{|\phi^{(\ell+1)}_{(\boldsymbol{k})}|^{p+1}
    - |\phi^{(\ell)}_{(\boldsymbol{k})}|^{p+1}}{
    \phi^{(\ell+1)}_{(\boldsymbol{k})}- \phi^{(\ell)}_{(\boldsymbol{k})}}
  \nonumber\\
  &\quad
  + \dfrac{e^{(n-2)Ht^{(\ell+1)}}+e^{(n-2)Ht^{(\ell)}}}{4}
  \delta^{ij}
  \widehat{\delta}^{\langle1\rangle}_i
  \widehat{\delta}^{\langle1\rangle}_j
  (\phi^{(\ell+1)}_{(\boldsymbol{k})} + \phi^{(\ell)}_{(\boldsymbol{k})}),
  \label{eq:DiscretePsi-I}
\end{align}
where ${}^{(\ell)}$ means the time index, ${}_{(\boldsymbol{k})}$ means the
space index, and $\boldsymbol{k}=(k_1,\dots,k_n)$.
$\widehat{\delta}^{\langle1\rangle}$ is the first-order central difference
operator defined as
\begin{align*}
  \widehat{\delta}^{\langle1\rangle}_{i}u^{(\ell)}_{(\boldsymbol{k})}
  := \dfrac{u^{(\ell)}_{(k_1,\dots,k_i+1,\dots,k_n)}
    - u^{(\ell)}_{(k_1,\dots,k_i-1,\dots,k_n)}}{2\Delta x^i}.
\end{align*}
Here, $\Delta x^i$ is the $i$-th grid range.
If $n=3$, for example, $\Delta x^1=\Delta x$, $\Delta x^2=\Delta y$, and
$\Delta x^3=\Delta z$.

The discretized \eqref{eq:totalH} can be defined as
\begin{align}
  {H}^{\mathrm{I}}_C{}^{(\ell)}
  :=
  \sum_{2\le k_1\le N_1+1}\cdots\sum_{2\le k_n\le N_n+1}
  \mathcal{H}^{\mathrm{I}}{}^{(\ell)}_{
    (\boldsymbol{k})}\Delta V,
  \label{eq:discreteTotalH-I}
\end{align}
where $\Delta V=\Delta x^1\cdots \Delta x^n$ and $N_i$ is the number of $i$-th
grids.
The discretized \eqref{eq:modTotalHam} is defined as
\begin{align}
  \tilde{H}^{\mathrm{I}}_C{}^{(\ell)}
  &:=
  H^{\mathrm{I}}_C{}^{(\ell)}-
  \sum_{q=0}^{\ell-1}(H^{\mathrm{I}}_C{}^{(q+1)} - H^{\mathrm{I}}_C{}^{(q)})
  \Delta t.
  \label{eq:disModTotalHamI}
\end{align}

\subsection{Form II}

The discretized \eqref{eq:HamiltonianDeSitter}, \eqref{eq:psi}, and
\eqref{eq:phi} can be respectively defined as
\begin{align}
  \mathcal{H}^{\mathrm{II}}{}^{(\ell)}_{(\boldsymbol{k})}
  &:= \frac{1}{2}e^{-nHt^{(\ell)}}
  (\psi^{(\ell)}_{(\boldsymbol{k})})^2
  + \dfrac{m^2}{2}e^{nHt^{(\ell)}}
  (\phi^{(\ell)}_{(\boldsymbol{k})})^2
  + \dfrac{1}{p+1}e^{nHt^{(\ell)}}
  |\phi^{(\ell)}_{(\boldsymbol{k})}|^{p+1}
  \nonumber\\
  &\quad
  + \dfrac{1}{4}e^{(n-2)Ht^{(\ell)}}\delta^{ij}\biggl(
  (\widehat{\delta}^{+}_i\phi^{(\ell)}_{(\boldsymbol{k})})
  (\widehat{\delta}^{+}_j\phi^{(\ell)}_{(\boldsymbol{k})})
  + (\widehat{\delta}^{-}_i\phi^{(\ell)}_{(\boldsymbol{k})})
  (\widehat{\delta}^{-}_j\phi^{(\ell)}_{(\boldsymbol{k})})
  \biggr),
  \label{eq:HamiltonianDeSitter-II}
\end{align}
\begin{align}
  \frac{\psi^{(\ell+1)}_{(\boldsymbol{k})}
    - \psi^{(\ell)}_{(\boldsymbol{k})}}{\Delta t}
  &:=
  - \dfrac{m^2}{4}(e^{nHt^{(\ell+1)}} + e^{nHt^{(\ell)}})
  (\phi^{(\ell+1)}_{(\boldsymbol{k})} + \phi^{(\ell)}_{(\boldsymbol{k})})
  - \frac{1}{2(p+1)}(e^{nHt^{(\ell+1)}}+ e^{nHt^{(\ell)}})
  \frac{|\phi^{(\ell+1)}_{(\boldsymbol{k})}|^{p+1}
    - |\phi^{(\ell)}_{(\boldsymbol{k})}|^{p+1}}{
    \phi^{(\ell+1)}_{(\boldsymbol{k})}
    - \phi^{(\ell)}_{(\boldsymbol{k})}}
  \nonumber\\
  &\quad
  + \dfrac{e^{(n-2)Ht^{(\ell+1)}}+e^{(n-2)Ht^{(\ell)}}}{4}
  \delta^{ij}\widehat{\delta}^{\langle2\rangle}_{ij}
  (\phi^{(\ell+1)}_{(\boldsymbol{k})} + \phi^{(\ell)}_{(\boldsymbol{k})}),
  \label{eq:DiscretePsi-II}
\end{align}
and \eqref{eq:DiscretePhi-I}, where $\widehat{\delta}^+_i$ is the first-order
forward difference operator defined as
\begin{align*}
  \widehat{\delta}^{+}_{i}u^{(\ell)}_{(\boldsymbol{k})}
  = \dfrac{u^{(\ell)}_{(k_1,\dots,k_i+1,\dots,k_n)}
    - u^{(\ell)}_{(\boldsymbol{k})}}{\Delta x^i},
\end{align*}
and $\widehat{\delta}^-_i$ is the first-order backward difference operator
defined as
\begin{align*}
  \widehat{\delta}^{-}_{i}u^{(\ell)}_{(\boldsymbol{k})}
  = \dfrac{u^{(\ell)}_{(\boldsymbol{k})}
    - u^{(\ell)}_{(k_1,\dots,k_i-1,\dots,k_n)}
    }{\Delta x^i}.
\end{align*}
$\widehat{\delta}^{\langle2\rangle}_{ij}$ is defined as
$\widehat{\delta}^{\langle2\rangle}_{ij}=\widehat{\delta}^{+}_{i}
\widehat{\delta}^{-}_{j}(=\widehat{\delta}^{-}_{j}
\widehat{\delta}^{+}_{i})$.
If $i=j$, the expression is the well-known second-order central difference
defined as
\begin{align*}
  \widehat{\delta}^{\langle2\rangle}_{ii}u^{(\ell)}_{(\boldsymbol{k})}
  = \dfrac{u^{(\ell)}_{(k_1\dots,k_i+1,\dots,k_n)}
    - 2u^{(\ell)}_{(\boldsymbol{k})}
    + u^{(\ell)}_{(k_1\dots,k_i-1,\dots,k_n)}}{(\Delta x^i)^2}.
\end{align*}

The discretized \eqref{eq:totalH} can be defined as
\begin{align}
  {H}^{\mathrm{II}}_C{}^{(\ell)}
  :=
  \sum_{1\le k_1\le N_1}\cdots
  \sum_{1\le k_n\le N_n}
  \mathcal{H}^{\mathrm{II}}{}^{(\ell)}_{
    (\boldsymbol{k})}\Delta V,
  \label{eq:discreteTotalH-II}
\end{align}
and the discretized \eqref{eq:modTotalHam} is defined as
\begin{align}
  \tilde{H}^{\mathrm{II}}_C{}^{(\ell)}
  &:=
  H^{\mathrm{II}}_C{}^{(\ell)}
  - \sum_{q=0}^{\ell-1}
  (H^{\mathrm{II}}_C{}^{(q+1)} - H^{\mathrm{II}}_C{}^{(q)})\Delta t.
  \label{eq:disModTotalHamII}
\end{align}
The values of \eqref{eq:discreteTotalH-I} and \eqref{eq:discreteTotalH-II}
are treated to investigate the numerical accuracy for $H=0$, and the values
of \eqref{eq:disModTotalHamI} and \eqref{eq:disModTotalHamII} are treated for
$H\neq 0$.

\section{Numerical simulations}

In this section, we perform some simulations with Forms I and II.
We set the initial conditions as  $\phi=A\cos(2\pi x)$ and $\psi=2\pi
A\sin(2\pi x)$, where $A=4$ and $-1/2\leq x\leq 1/2$.
The boundary is periodic.
The grid range is $(\Delta x, \Delta t)=(1/50, 1/250),(1/100, 1/500),(1/200,
1/1000)$.
The number of exponents of the nonlinear term is $p=5$.

\subsection{Flat spacetime}

\begin{figure}[t]  \centering
  \begin{minipage}{0.49\hsize}
  \includegraphics[width=1\hsize]{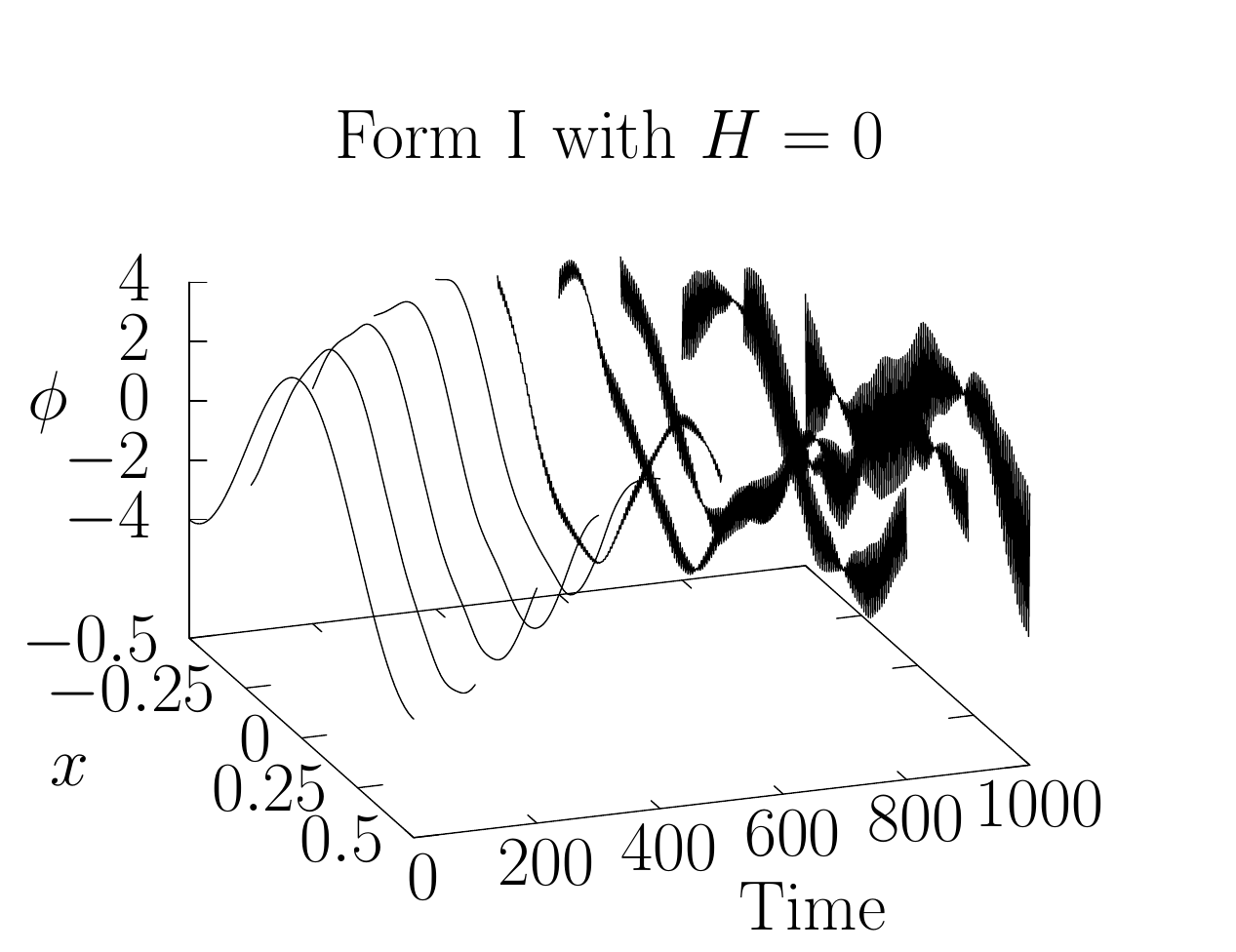}
  \end{minipage}
  \begin{minipage}{0.49\hsize}
  \includegraphics[width=1\hsize]{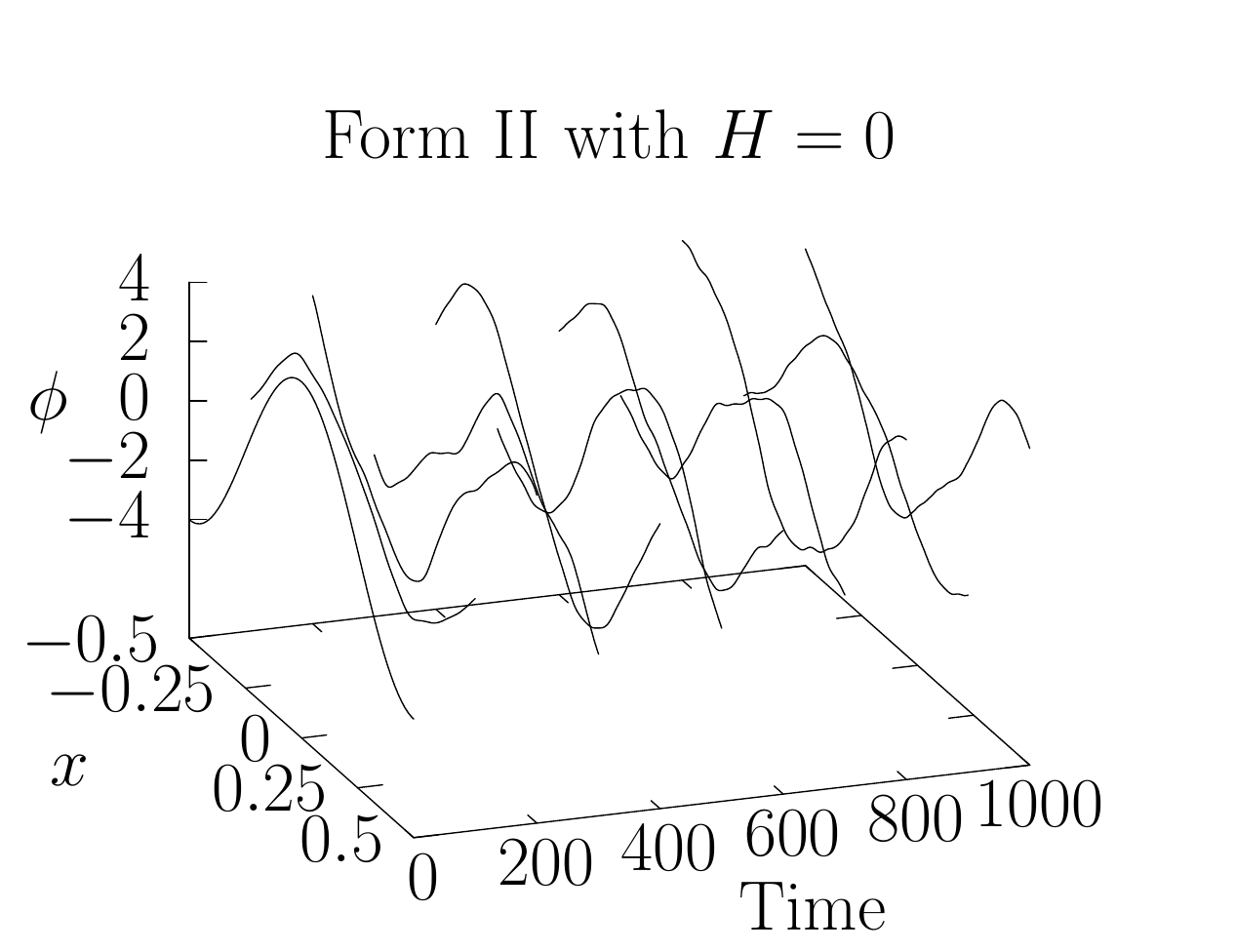}
  \end{minipage}
  \caption{%
    \label{fig:H0phi}
    $\phi$ obtained with Forms I and II in $H=0$.
    The left panel is drawn with Form I and the right one with Form II.
    The vibrations occur after $t\geq 500$ in the left panel.}
\end{figure}

First, we perform simulations in a flat spacetime, that is, $H=0$.
Fig. \ref{fig:H0phi} shows the waveform of $\phi$ obtained with Forms I and II.
We see that the vibrations occur in the left panel, which is drawn with Form I,
and no vibrations occur in the right panel, which is drawn with Form II.
This means that the simulation with Form II is more stable than that with Form
I.

\begin{figure}[t]
  \centering
  \includegraphics[width=0.65\hsize]{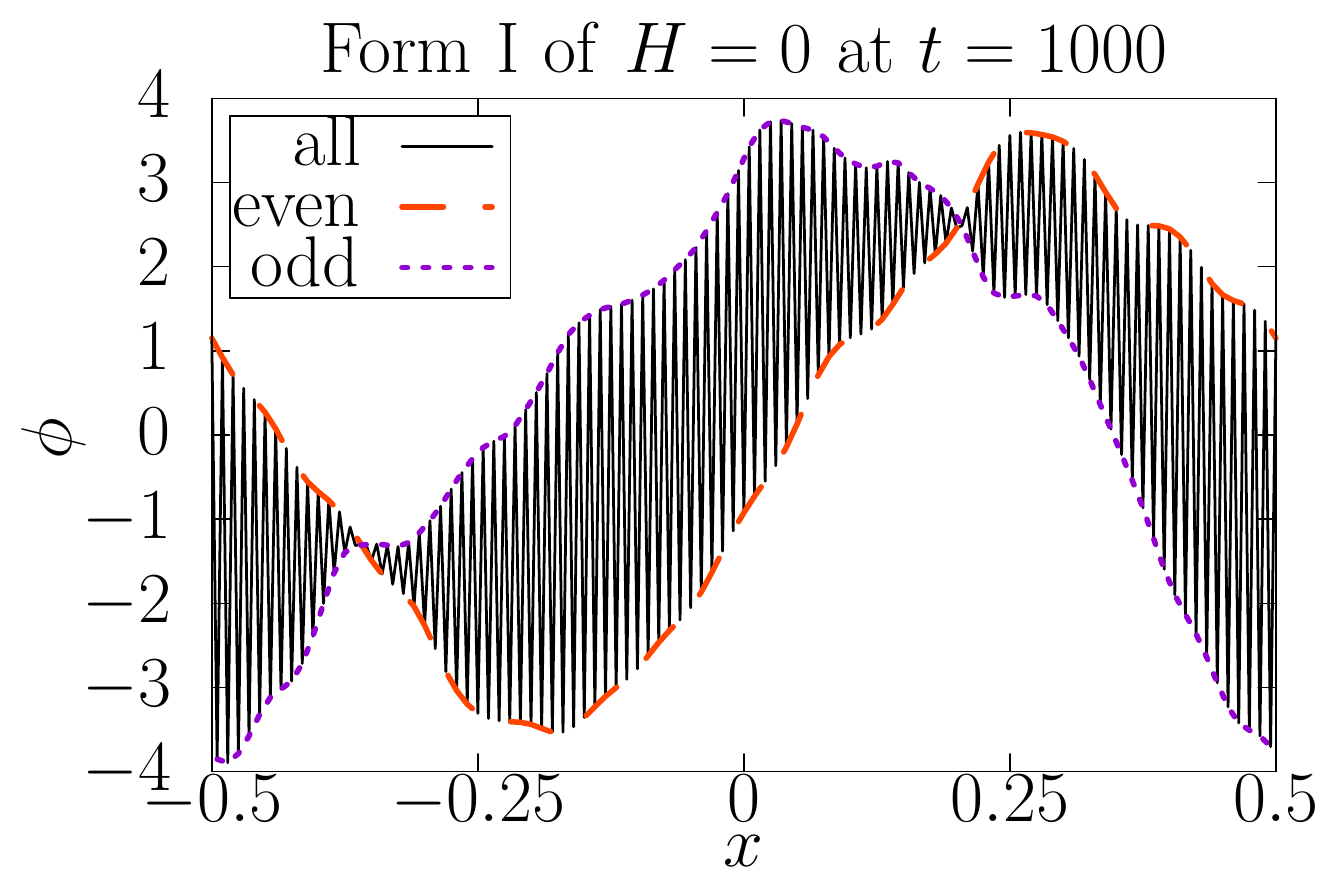}
  \caption{%
    \label{fig:H0phi1000}
    Waveform of $\phi$ obtained with Form I at $t=1000$ for $\Delta x=1/200$ and
    $\Delta t=1/1000$.
    The dashed and dotted lines are drawn with the even and odd numbers of grid
    points, respectively.
    On the other hand, the solid line is drawn with all grid points.}
\end{figure}

We determine the reasons for the generation of vibrations in the waveform with
Form I.
Fig. \ref{fig:H0phi1000} shows the waveform of $\phi$ obtained with Form I at
$t=1000$.
There are marked differences between the line with an even number of grid points
and that with an odd number of grid points.
These are mainly caused by the second-order difference formula in
\eqref{eq:DiscretePsi-I}.
This is explicitly expressed as
\begin{align*}
  \delta^{ij}\widehat{\delta}^{\langle1\rangle}_i
  \widehat{\delta}^{\langle1\rangle}_j\phi^{(n)}_{(\boldsymbol{k})}
  &= \sum_{i=1}^n\dfrac{\phi^{(n)}_{(k_1,\dots,k_i+2,\dots,k_n)}
    - 2\phi^{(n)}_{(\boldsymbol{k})}
    + \phi^{(n)}_{(k_1,\dots,k_i-2,\dots,k_n)}}{(\Delta x^i)^2}.
\end{align*}
The expression indicates that the odd and even numbers of grid points are
independent of each other.
If the differences between the odd- and even-number grid points are generated by
numerical errors, the grid points separate into those with even and odd numbers.
On the other hand, the formula in \eqref{eq:DiscretePsi-II} is explicitly
expressed as
\begin{align*}
  \delta^{ij}\widehat{\delta}^{\langle2\rangle}_{ij}
  \phi^{(n)}_{(\boldsymbol{k})}
  &= \sum_{i=1}^n\dfrac{\phi^{(n)}_{(k_1,\dots,k_i+1,\dots,k_n)}
    - 2\phi^{(n)}_{(\boldsymbol{k})}
    + \phi^{(n)}_{(k_1,\dots,k_i-1,\dots,k_n)}}{(\Delta x^i)^2}.
\end{align*}
Thus, the odd- and even-number grid points are dependent on each other.
If the differences are generated by numerical errors, the differences propagate
at all grid points.
Therefore, no vibrations occur in $\phi$ obtained with Form II as shown in Fig.
\ref{fig:H0phi}.

\begin{figure}[t]
  \centering
  \begin{minipage}{0.49\hsize}
  \includegraphics[width=\hsize]{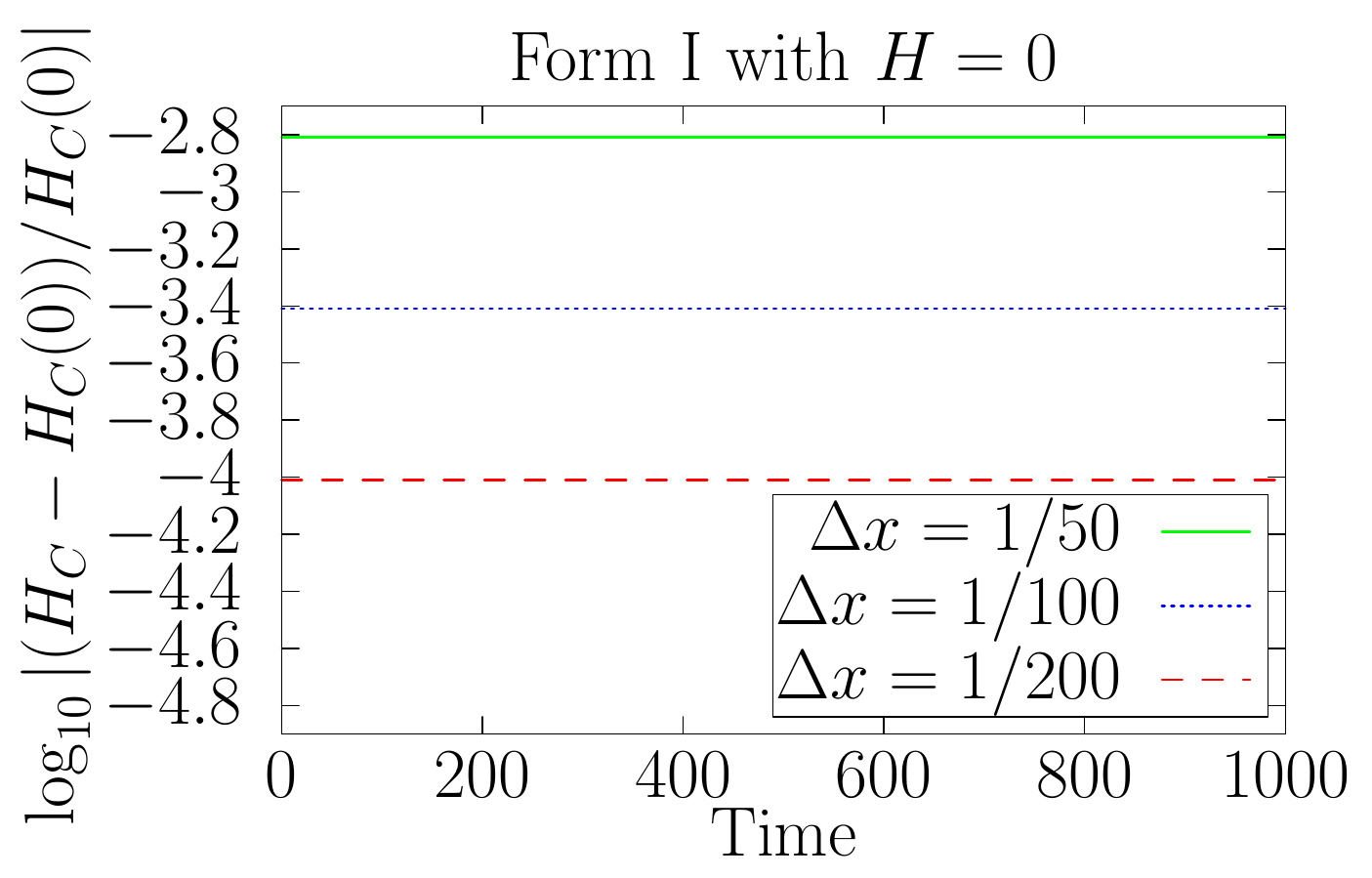}
  \end{minipage}
  \begin{minipage}{0.49\hsize}
    \includegraphics[width=\hsize]{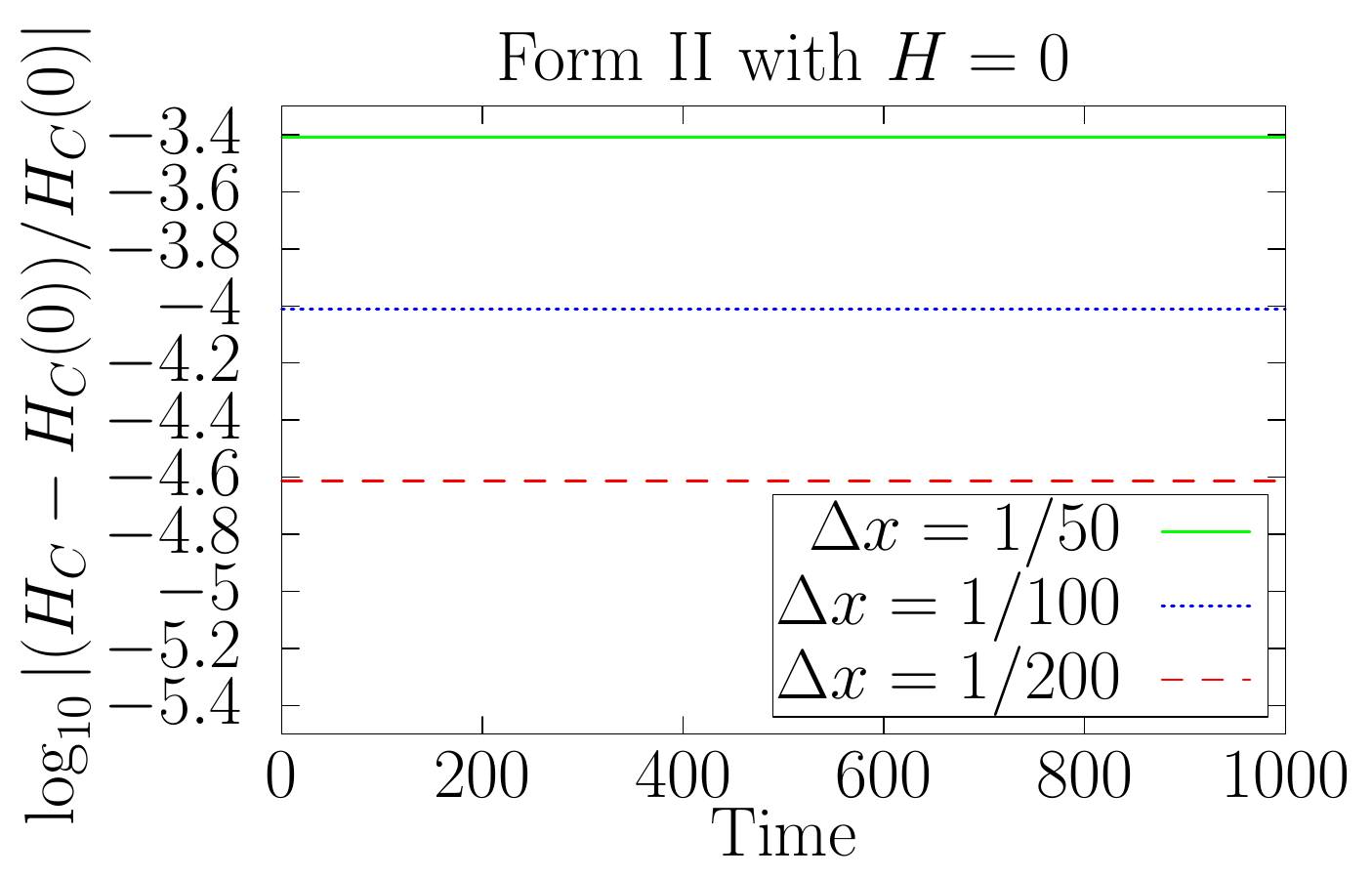}
  \end{minipage}
  \caption{
    \label{fig:HC}
    Relative errors of $H^{\mathrm{I}}_C$ against the initial value obtained
    with Form I and those of $H^{\mathrm{II}}_C$ against the initial value
    obtained with Form II in the Hubble constant as $H=0$.
    The left panel is drawn with Form I and the right one with Form II.}
\end{figure}
Fig. \ref{fig:HC} shows relative errors of the discretized total Hamiltonians
$H^{\mathrm{I}}_C$ and $H^{\mathrm{II}}_C$ against the initial value with Forms
I and II, respectively.
When the number of grid points is increased twofold, the value is about
$0.6\approx \log_{10}4$ smaller in both panels.
The results mean that $H^{\mathrm{I}}_C$ and $H^{\mathrm{II}}_C$ show the
second-order accuracies with respect to the number of grid points.
However, $H^{\mathrm{II}}_C$ should show the first-order accuracy because of the
expression of \eqref{eq:HamiltonianDeSitter-II}.
This discrepancy will be discussed in Sec. \ref{sec:sum_dis}.

\subsection{Curved spacetime}
\begin{figure}[t]
  \centering
  \begin{minipage}{0.49\hsize}
  \includegraphics[width=\hsize]{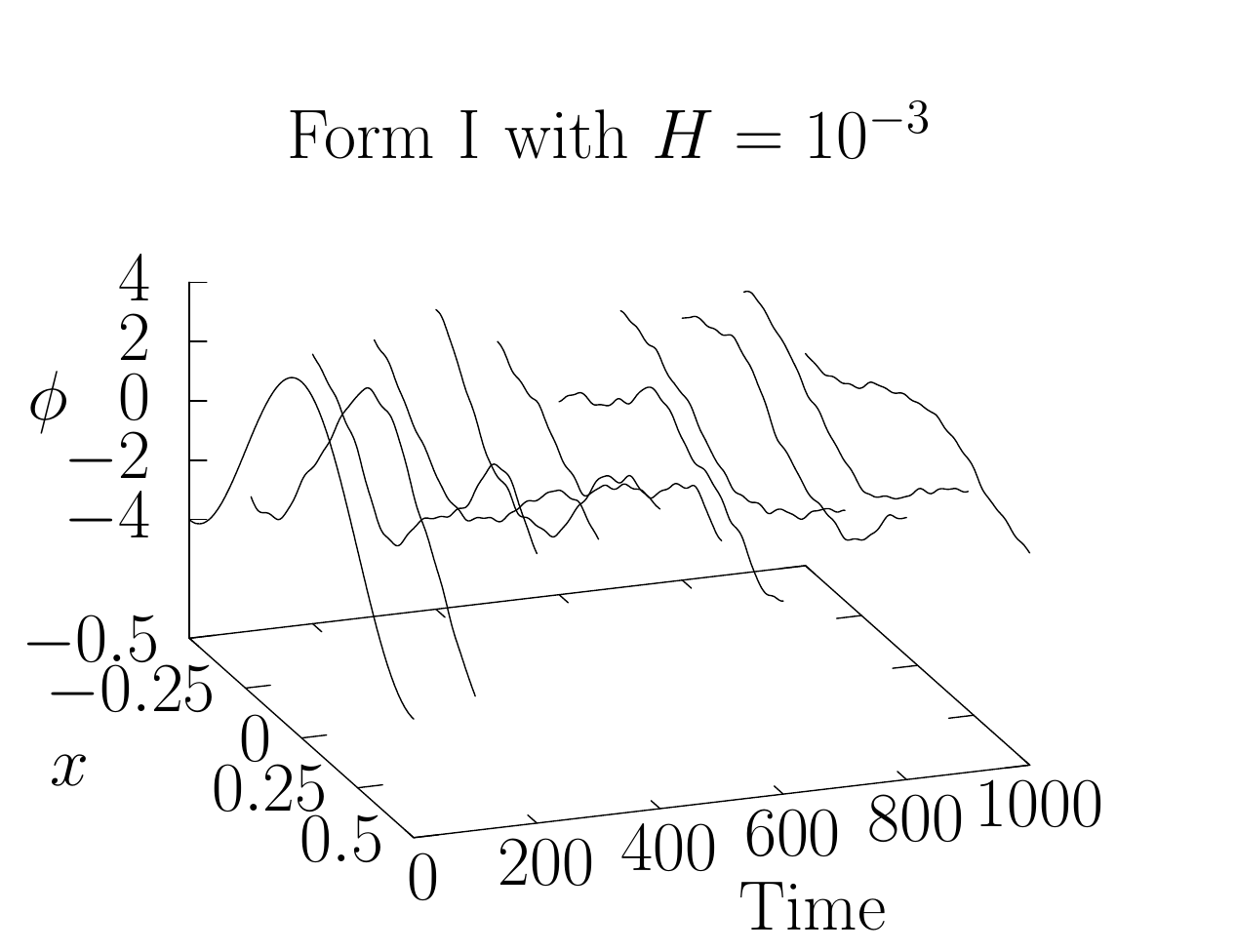}
  \end{minipage}
  \begin{minipage}{0.49\hsize}
  \includegraphics[width=\hsize]{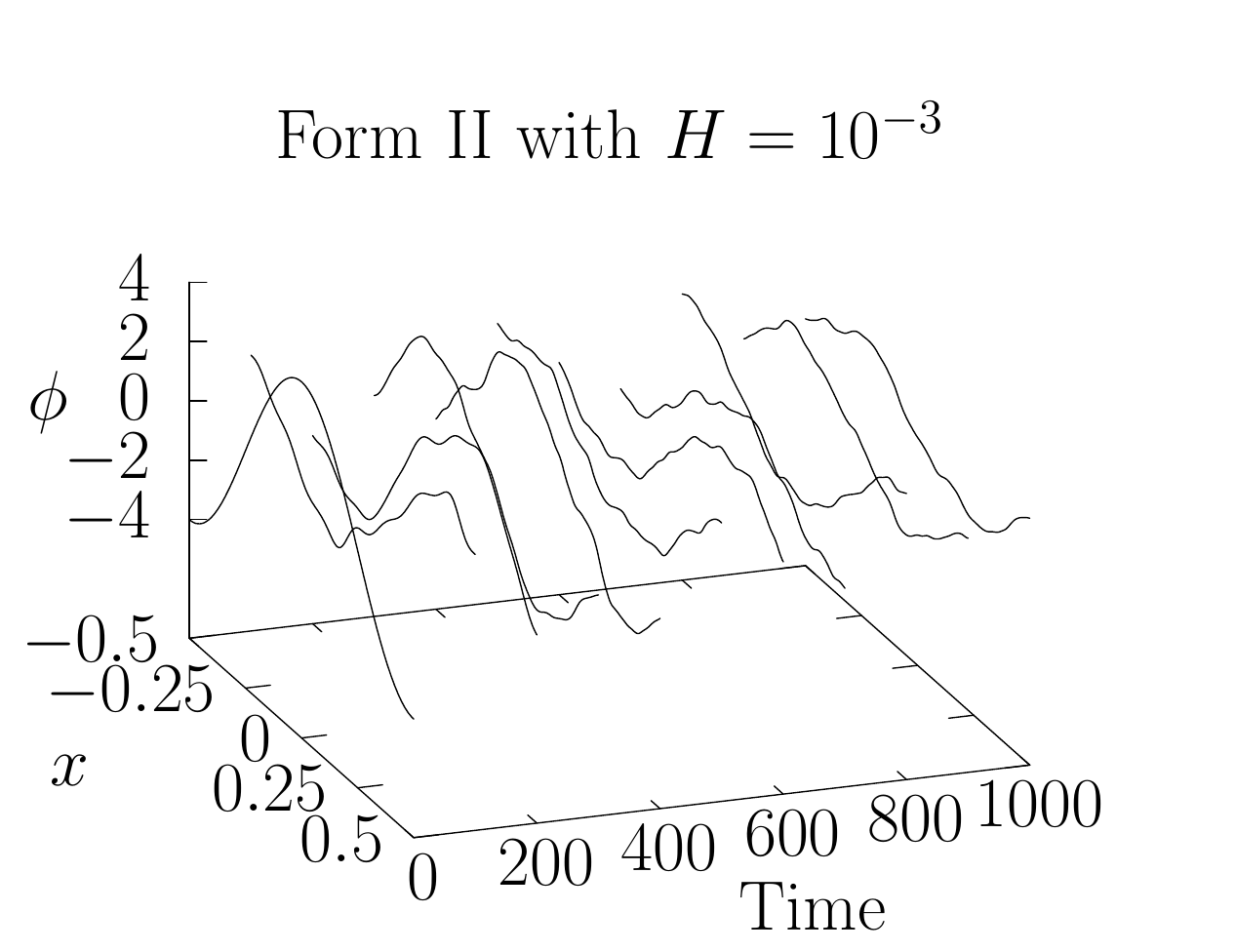}
  \end{minipage}
  \caption{
    \label{fig:CurvedPhi}
    $\phi$ obtained with Forms I and II in $H=10^{-3}$.
    The left panel is drawn with Form I and the right one with Form II.}
\end{figure}

\begin{figure}[t]
  \centering
  \begin{minipage}{0.49\hsize}
  \includegraphics[width=\hsize]{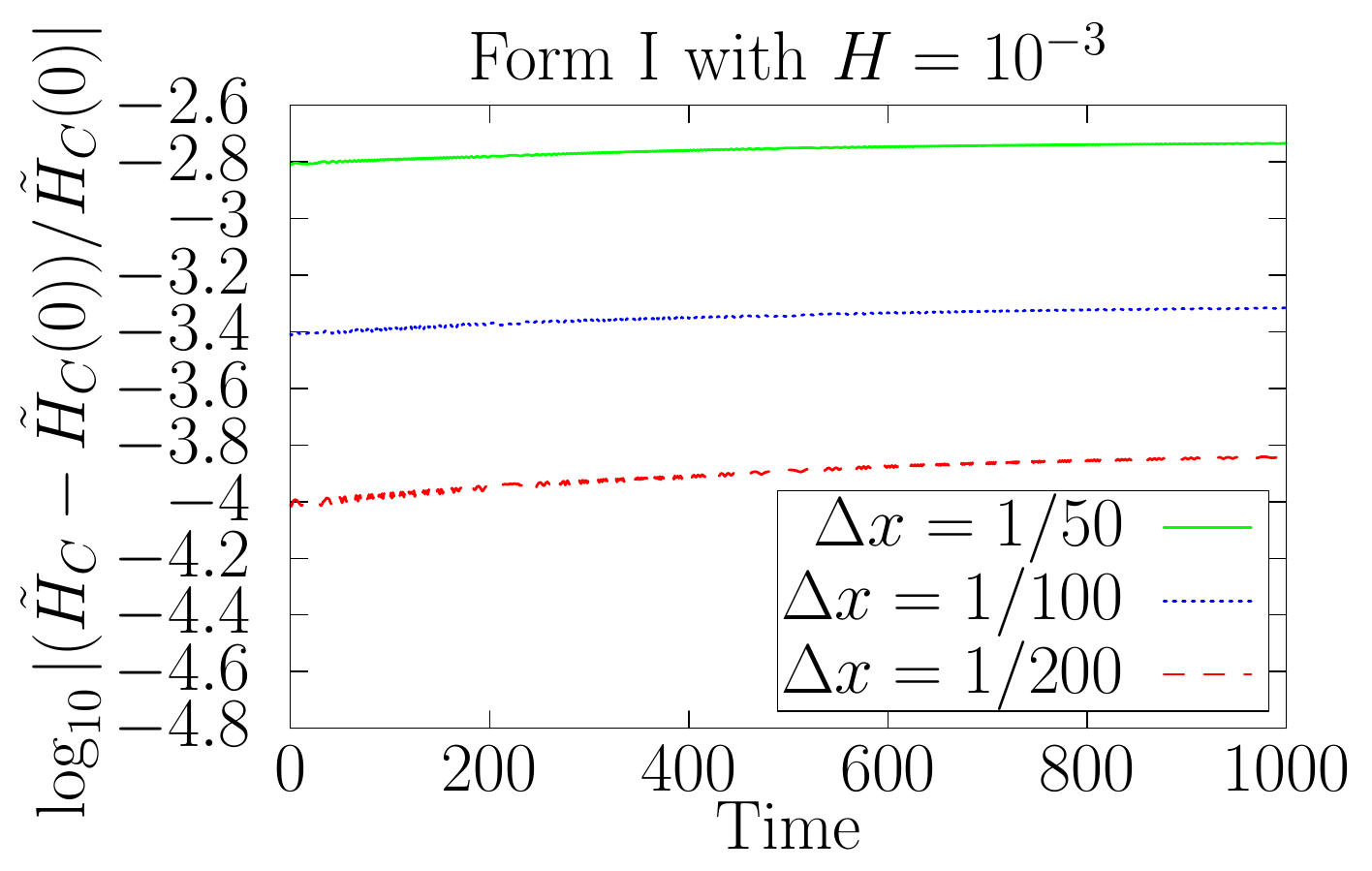}
  \end{minipage}
  \begin{minipage}{0.49\hsize}
  \includegraphics[width=\hsize]{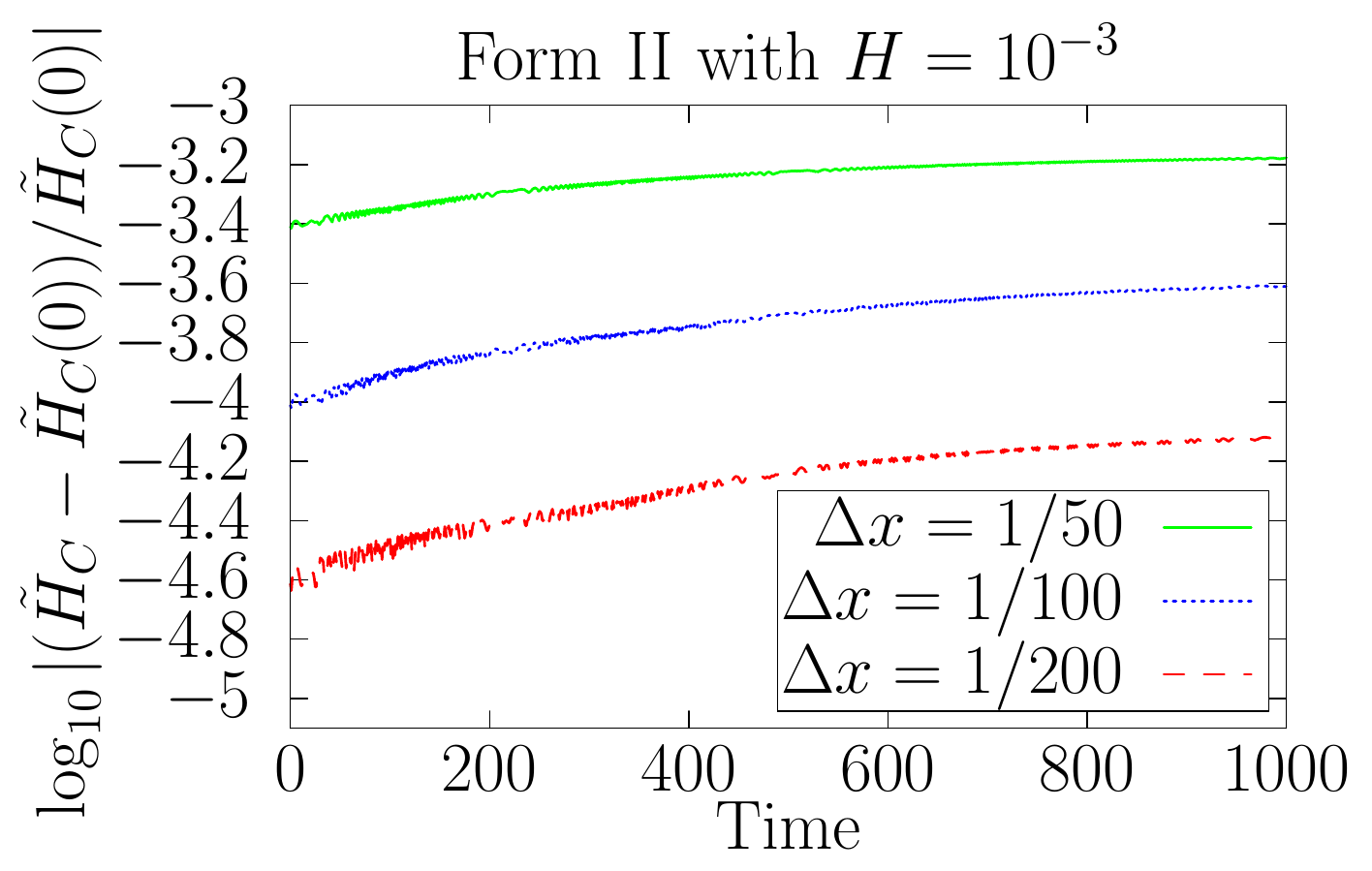}
  \end{minipage}
  \caption{
    \label{fig:CurvedHC}
    Relative errors of $\tilde{H}^{\mathrm{I}}_C$ against the initial value
    obtained with Form I and those of $\tilde{H}^{\mathrm{II}}_C$ against the
    initial value obtained with Form II in the Hubble constant as $H=10^{-3}$.
    The left panel is drawn with Form I and the right one with Form II.}
\end{figure}

Next, we perform some simulations in an expanding space, that is, $H=10^{-3}$.
The purpose of this study is to confirm the efficiency of the Hubble constant
$H$ in terms of the accuracy and stability of the simulations.
In Fig. \ref{fig:CurvedPhi}, no vibrations appear.
This has already been mentioned in \cite{Tsuchiya-Nakamura-2019-JCAM}, which
means that the expansions of the space increase the stability of the
simulations.
The modified total Hamiltonian $\tilde{H}_C$ is drawn in Fig.
\ref{fig:CurvedHC}.
At the initial time, the numerical accuracy is almost of the second order with
respect to the grid.
However, near $t=1000$, the accuracy is less than the second order.
Therefore, the numerical accuracy decreases with time if $H\neq 0$.

%
\section{Summary and discussion
  \label{sec:sum_dis}
}

We performed some simulations of the semilinear Klein--Gordon equation in the
de Sitter spacetime with two structure-preserving discretized forms of the
equation.
Form I has the product of the first-order central difference formulae in the
discretized evolution equations.
Form II has the product of the first-order forward and backward difference
formulae.
If $H=0$, Form II is more stable than Form I because Form I but not Form II has
numerical vibrations in the waveform.
On the other hand, if $H=10^{-3}$, there are no vibrations in the waveform.

The numerical accuracies of the two forms are of the second order with respect
to the grid if $H=0$.
However, the expression of the total Hamiltonian of Form I indicates
first-order accuracy.

In the construction from \eqref{eq:HamiltonianDeSitter} to
\eqref{eq:HamiltonianDeSitter-II}, the differential term in the one-dimensional
case is
\begin{align*}
  (\partial_x\phi)(\partial_x\phi)
  &=
  \biggl(\dfrac{\phi^{(n)}_{(k+1)} - \phi^{(n)}_{(k)}}{\Delta x}
  - \dfrac{\Delta x}{2}  \phi''+O((\Delta x)^2)
  \biggr)^2
  + \biggl(\dfrac{\phi^{(n)}_{(k)} - \phi^{(n)}_{(k-1)}}{\Delta x}
  + \dfrac{\Delta x}{2} \phi''+O((\Delta x)^2)
  \biggr)^2
  \nonumber\\
  &= \biggl(\dfrac{\phi^{(n)}_{(k+1)} - \phi^{(n)}_{(k)}}{\Delta x}\biggr)^2
  + \biggl(\dfrac{\phi^{(n)}_{(k)} - \phi^{(n)}_{(k-1)}}{\Delta x}\biggr)^2
  - \phi''(\phi^{(n)}_{(k+1)} - 2\phi^{(n)}_{(k)}
    + \phi^{(n)}_{(k-1)})
  +O((\Delta x)^2)
  \nonumber\\
  &=
  \biggl(\dfrac{\phi^{(n)}_{(k+1)} - \phi^{(n)}_{(k)}}{\Delta x}\biggr)^2
  + \biggl(\dfrac{\phi^{(n)}_{(k)} - \phi^{(n)}_{(k-1)}}{\Delta x}\biggr)^2
  +O((\Delta x)^2),
\end{align*}
where we use the relation given by the evolution equation
\eqref{eq:DiscretePsi-II} such as
\begin{align*}
  \phi^{(n)}_{(k+1)} - 2\phi^{(n)}_{(k)} + \phi^{(n)}_{(k-1)}
  &= (\Delta x)^2\biggl[
    \dfrac{\psi^{(n+1)}_{(k)}-\psi^{(n)}_{(k)}}{\Delta t}
    + m^2\phi^{(n)}_{(k)}
    + |\phi^{(n)}_{(k)}|^{p-1}
    \biggr]
  \\
  &= O((\Delta x)^2).
\end{align*}
Therefore, we conclude that the relative errors of $H^{\mathrm{II}}_C$ obtained
with Form II might affect the second-order accuracy with respect to the grid.

\section*{Acknowledgments}
T.T. and M.N. were partially supported by JSPS KAKENHI Grant Number 21K03354.
T.T. was partially supported by JSPS KAKENHI Grant Number 20K03740
and Grant for Basic Science Research Projects from The Sumitomo Foundation.
M.N. was partially supported by JSPS KAKENHI Grant Number 16H03940.


\end{document}